# $L^P$ MODULI OF CONTINUITY OF GAUSSIAN PROCESSES AND LOCAL TIMES OF SYMMETRIC LÉVY PROCESSES[1]

BY MICHAEL B. MARCUS AND JAY ROSEN

*City University of New York*

Let $X = \{X(t), t \in R_+\}$ be a real-valued symmetric Lévy process with continuous local times $\{L_t^x, (t,x) \in R_+ \times R\}$ and characteristic function $Ee^{i\lambda X(t)} = e^{-t\psi(\lambda)}$. Let

$$\sigma_0^2(x-y) = \frac{4}{\pi} \int_0^\infty \frac{\sin^2(\lambda(x-y)/2)}{\psi(\lambda)} \, d\lambda.$$

If $\sigma_0^2(h)$ is concave, and satisfies some additional very weak regularity conditions, then for any $p \geq 1$, and all $t \in R_+$,

$$\lim_{h \downarrow 0} \int_a^b \left| \frac{L_t^{x+h} - L_t^x}{\sigma_0(h)} \right|^p dx = 2^{p/2} E|\eta|^p \int_a^b |L_t^x|^{p/2} \, dx$$

for all $a, b$ in the extended real line almost surely, and also in $L^m$, $m \geq 1$. (Here $\eta$ is a normal random variable with mean zero and variance one.)

This result is obtained via the Eisenbaum Isomorphism Theorem and depends on the related result for Gaussian processes with stationary increments, $\{G(x), x \in R^1\}$, for which $E(G(x) - G(y))^2 = \sigma_0^2(x-y)$;

$$\lim_{h \to 0} \int_a^b \left| \frac{G(x+h) - G(x)}{\sigma_0(h)} \right|^p dx = E|\eta|^p (b-a)$$

for all $a, b \in R^1$, almost surely.

**1. Introduction.** We obtain $L^p$ moduli of continuity for a very wide class of continuous Gaussian processes and local times of symmetric Lévy processes. To introduce them, we first state our results for the local times of the Brownian motion and see how they compare with related results.

Received July 2006; revised December 2006.

[1]Supported in part by grants from the National Science Foundation and PSC-CUNY.

THEOREM 1.1. *Let $\{L_t^x, (x,t) \in R^1 \times R_+\}$ denote the local time of Brownian motion. Then, for any $p \geq 1$ and $t \in R_+$,*

$$(1.1) \qquad \lim_{h \downarrow 0} \int_a^b \left| \frac{L_t^{x+h} - L_t^x}{h^{1/2}} \right|^p dx = \frac{2^{3p/2}}{\sqrt{\pi}} \Gamma\left(\frac{p+1}{2}\right) \int_a^b |L_t^x|^{p/2} dx$$

*for all $a, b$ in the extended real line $y$.*

When $p = 2$, (1.1) is the following: For all $t \in R_+$,

$$(1.2) \qquad \lim_{h \downarrow 0} \frac{\int_{-\infty}^\infty (L_t^{x+h} - L_t^x)^2 dx}{h} = 4t \qquad \text{a.s.}$$

This may be considered as a continuous version of the quadratic variation result: For all $t \in R_+$,

$$(1.3) \qquad \lim_{n \to \infty} \sum_{j=-\infty}^\infty (L_t^{j/n} - L_t^{(j-1)/n})^2 = 4 \int_{-\infty}^\infty L_t^x \, dx = 4t \qquad \text{a.s.}$$

(We obtain (1.3) from [2], Theorem 10.4.1 and Lemma 10.5.2, using the 2-stable process which is the Brownian motion multiplied by $\sqrt{2}$.)

When $p = 1$, (1.1) is the following: For all $t \in R_+$,

$$(1.4) \qquad \lim_{h \downarrow 0} \frac{\int_a^b |L_t^{x+h} - L_t^x| \, dx}{\sqrt{h}} = \frac{2^{3/2}}{\sqrt{\pi}} \int_a^b \sqrt{L_t^x} \, dx \qquad \text{a.s.}$$

This compliments a result of Yor [4] that

$$(1.5) \qquad \lim_{h \downarrow 0} \frac{L_t^h - L_t^0}{\sqrt{h}} \stackrel{\text{law}}{=} 2\sqrt{L_t^0} \eta,$$

where $\eta$ is a normal random variable with mean zero and variance one.

Theorem 1.1 can be extended to symmetric Lévy processes with continuous local times, subject to some regularity conditions. Let $X = \{X(t), t \in R_+\}$ be a real-valued symmetric Lévy process with characteristic function

$$(1.6) \qquad E e^{i\lambda X(t)} = e^{-t\psi(\lambda)},$$

where

$$(1.7) \qquad \psi(\lambda) = 2 \int_0^\infty (1 - \cos \lambda u) \nu(du)$$

for $\nu$ a symmetric Lévy measure, that is, $\nu$ is symmetric and

$$(1.8) \qquad \int_0^\infty (1 \wedge x^2) \nu(dx) < \infty.$$

We assume that

$$(1.9) \qquad \int_1^\infty \frac{1}{\psi(\lambda)} d\lambda < \infty,$$

which is a necessary and sufficient condition for $X$ to have local times. We refer to $\psi(\lambda)$ as the characteristic exponent of $X$. Let

$$\sigma_0^2(x-y) = \frac{4}{\pi} \int_0^\infty \frac{\sin^2(\lambda(x-y)/2)}{\psi(\lambda)} d\lambda. \tag{1.10}$$

We say that $\sigma_0$ satisfies condition $\mathcal{C}_q$ if

$$\lim_{n\to\infty} \frac{\sigma_0(1/n(\log n)^{q+1})}{\sigma_0(1/(\log n)^q)} = 0. \tag{1.11}$$

We say that $\psi(\lambda)$ satisfies condition $\Lambda_\gamma$ if

$$\lambda^\gamma = o(\psi(\lambda)) \qquad \text{as } \lambda \to \infty. \tag{1.12}$$

THEOREM 1.2. *Let $X = \{X(t), t \in R_+\}$ be a real-valued symmetric Lévy process with characteristic exponent $\psi(\lambda)$ that satisfies condition $\Lambda_\gamma$, for some $\gamma > 0$. Assume that $\sigma_0^2(h)$ is concave and monotonically increasing for $h \in [0, \delta]$ for some $\delta > 0$ and satisfies condition $\mathcal{C}_q$. Let $L := \{L_t^x, (t,x) \in R_+ \times R\}$ be the local time of $X$ and assume that $L$ is continuous. Let $\eta$ be a normal random variable with mean zero and variance one. Then for any $1 \le p < q$ and all $t \in R_+$,*

$$\lim_{h\downarrow 0} \int_a^b \left|\frac{L_t^{x+h} - L_t^x}{\sigma_0(h)}\right|^p dx = 2^{p/2} E|\eta|^p \int_a^b |L_t^x|^{p/2} dx$$
$$= \frac{2^p}{\sqrt{\pi}} \Gamma\left(\frac{p+1}{2}\right) \int_a^b |L_t^x|^{p/2} dx \tag{1.13}$$

*for all $a, b$ in the extended real line almost surely.*

We point out on page 22 for which $\sigma_0^2$ is concave. The other two conditions in this theorem are very weak.

In Section 5 we show that the limit in (1.13) also exists in $L^m$ uniformly in $t$ on any bounded interval of $R_+$, for all $m \ge 1$.

When $\psi(\lambda) = |\lambda|^\beta$, $1 < \beta \le 2$, we refer to $X$ as the canonical $\beta$-stable process. (The canonical 2-stable process is the Brownian motion multiplied by $\sqrt{2}$.) In this case the conditions in Theorem 1.2 hold and (1.13) is the following: For any $1 \le p < q$ and $t \in R_+$,

$$\lim_{h\downarrow 0} \int_a^b \frac{|L_t^{x+h} - L_t^x|^p}{h^{p(\beta-1)/2}} dx = c(\beta, p) \int_a^b |L_t^x|^{p/2} dx \tag{1.14}$$

for all $a, b$ in the extended real line almost surely, where

$$c(\beta, p) = \left(\frac{1}{\Gamma(\beta)\sin((\pi/2)(\beta-1))}\right)^{p/2} \frac{2^p}{\sqrt{\pi}} \Gamma\left(\frac{p+1}{2}\right). \tag{1.15}$$

(See Remark 4.1 for more details.)

We derive our results on the $L^p$ moduli of continuity of local times of symmetric Lévy processes using the Eisenbaum Isomorphism Theorem ([2], Theorem 8.1.1). In order to use it, we need to know about the $L^p$ moduli of continuity of squares of the associated Gaussian processes. These follow easily from results about the $L^p$ moduli of continuity of the Gaussian processes themselves. These are interesting in their own right. We take this up in the next section. Here we just mention an application of the results to the fractional Brownian motion. Let $G = \{G(x), x \in R^1\}$ be a real-valued Gaussian process with mean zero and stationary increments, $G(0) = 0$, and let

$$E(G(x+h) - G(x))^2 = h^r, \tag{1.16}$$

$0 < r < 2$. Then

$$\lim_{h \downarrow 0} \int_a^b \left| \frac{G(x+h) - G(x)}{h^{r/2}} \right|^p dx = E|\eta|^p (b-a) \tag{1.17}$$

for all $-\infty < a < b < \infty$ almost surely. Results like (1.17) also follow from the work of Wschebor [3]. We explain in Remark 2.1 why we cannot use his approach to obtain Theorem 1.2.

**2. $L^p$ moduli of continuity of Gaussian processes.** Let $G = \{G(x), x \in R^1\}$ be a real-valued Gaussian process with mean zero and stationary increments and let

$$\sigma^2(h) = E(G(x+h) - G(x))^2. \tag{2.1}$$

Fix $1 \leq p < \infty$, $-\infty < a < b < \infty$ and define

$$I(h) = I_G(h; a, b, p) = \int_a^b \left| \frac{G(x+h) - G(x)}{\sigma(h)} \right|^p dx. \tag{2.2}$$

Then, clearly,

$$EI_G(h; a, b, p) = E|\eta|^p (b-a), \tag{2.3}$$

where $\eta$ is a normal random variable with mean zero and variance one. This shows, in particular, that $I_G(h; a, b, p)$ exists and is finite for all measurable Gaussian processes $G$. When $\sigma^2$ is concave in some neighborhood of the origin, $I_G(h; a, b, p)$ exhibits the following remarkable regularity property, whether $G$ has continuous paths or is unbounded almost surely. (These are the only two possibilities for $G$; see, e.g., [2], Theorem 5.3.10.)

THEOREM 2.1. *Let $G$ be as above and assume that $\sigma^2(h)$ is concave and monotonically increasing for $h \in [0, \delta]$, for some $\delta > 0$. Let $\{h_n\}$ be positive numbers with $h_n = o(\frac{1}{(\log n)^p})$. Then for any $1 \leq p < \infty$,*

$$(2.4) \qquad \lim_{n \to \infty} \int_a^b \left| \frac{G(x + h_n) - G(x)}{\sigma(h_n)} \right|^p dx = E|\eta|^p (b - a)$$

*for all $a, b \in R^1$, almost surely.*

Before proving this theorem, we give a preliminary lemma that is an application of the Borell, Sudakov–Tsirelson Theorem. For each $h$, consider the symmetric positive definite kernel

$$(2.5) \qquad \rho_h(x, y) = \frac{1}{\sigma^2(h)} E(G(x + h) - G(x))(G(y + h) - G(y)),$$

$$x, y \in R^1.$$

Note that by stationarity and the Cauchy–Schwarz inequality

$$(2.6) \qquad |\rho_h(x, y)| \leq 1, \qquad x, y \in R^1.$$

For $p \geq 1$, define

$$(2.7) \qquad \|G\|_{h,p} = (I_G(h; a, b, p))^{1/p}.$$

We denote the median of a real-valued random variable, say, $Z$, by $\mathrm{med}(Z)$.

LEMMA 2.1. *Under the hypotheses of Theorem 2.1,*

$$(2.8) \qquad P(|\|G\|_{h,p} - \mathrm{med}(\|G\|_{h,p})| > t) \leq 2 e^{-t^2/(2\widehat{\sigma}^2)},$$

*where*

$$(2.9) \qquad \widehat{\sigma}^2 = \sup_{\{f : \|f\|_q \leq 1\}} \int_a^b \int_a^b f(x) f(y) \rho_h(x, y) \, dx \, dy$$

*and $1/p + 1/q = 1$. Furthermore,*

$$(2.10) \qquad \widehat{\sigma}^2 \leq \left( \int_a^b \int_a^b |\rho_h(x, y)| \, dx \, dy \right)^{1/p}$$

*and*

$$(2.11) \qquad |E(\|G\|_{h,p}) - \mathrm{med}(\|G\|_{h,p})| \leq \frac{\widehat{\sigma}}{\sqrt{2\pi}}.$$

PROOF. Let $B_q$ be a countable dense subset of the unit ball of $L^q([a,b])$. For $f \in B_q$, set

$$H(h, f) = \int_a^b f(x) \frac{(G(x+h) - G(x))}{\sigma(h)} \, dx. \tag{2.12}$$

It is a standard fact in Banach space theory that

$$\sup_{f \in B_q} H(h, f) = \|G\|_{h,p}. \tag{2.13}$$

Let

$$\begin{aligned}
\widehat{\sigma}^2 &:= \sup_{f \in B_q} E(H^2(h, f)) \\
&= \sup_{\{f : \|f\|_q \leq 1\}} \int_a^b \int_a^b f(x) f(y) \rho_h(x, y) \, dx \, dy.
\end{aligned} \tag{2.14}$$

The statements in (2.8) and (2.9) follow from a standard application of the Borell, Sudakov–Tsirelson Theorem (see [2], Theorem 5.4.3).

For $1 \leq p < \infty$,

$$\begin{aligned}
\widehat{\sigma}^2 &\leq \left( \int_a^b \int_a^b |\rho_h(x, y)|^p \, dx \, dy \right)^{1/p} \\
&\leq \left( \int_a^b \int_a^b |\rho_h(x, y)| \, dx \, dy \right)^{1/p},
\end{aligned} \tag{2.15}$$

where in the last line we use (2.6). This follows from Hölder's inequality when $1 < p < \infty$. When $p = 1$, $q = \infty$ and $\|f\|_\infty := \sup_x |f(x)|$. Obtaining (2.15) in this case is trivial.

The statement in (2.11) is another standard application of the Borell, Sudakov–Tsirelson Theorem (see [2], Corollary 5.4.5). □

PROOF OF THEOREM 2.1. In order to use the concavity of $\sigma^2(h)$ on $[0, \delta]$, we initially take $b - a < \delta/2$. It follows from (2.8) and (2.10) that

$$P(|\|G\|_{h_n,p} - \mathrm{med}(\|G\|_{h_n,p})| > t) \leq 2 e^{-t^2/(2\widehat{\sigma}_n^2)}, \tag{2.16}$$

where

$$\widehat{\sigma}_n^2 \leq \left( \int_a^b \int_a^b |\rho_{h_n}(x, y)| \, dx \, dy \right)^{1/p}. \tag{2.17}$$

We show below that

$$\int_a^b \int_a^b |\rho_{h_n}(x, y)| \, dx \, dy = o\left( \frac{1}{(\log n)^p} \right) \tag{2.18}$$

as $n \to \infty$. Assuming this, we see from (2.16), (2.17), (2.18) and the Borel–Cantelli Lemma that

$$\lim_{n \to \infty} (\|G\|_{h_n,\rho} - \operatorname{med}\|G\|_{h_n,\rho}) = 0 \quad \text{a.s.} \tag{2.19}$$

Let $\operatorname{med}(\|G\|_{h_n,p}) = M_n$ and note that by (2.3)

$$M_n \leq 2E(\|G\|_{h_n,p}) \leq 2(E\|G\|_{h_n,p}^p)^{1/p}$$
$$= 2(E|\eta|^p)^{1/p}(b-a)^{1/p} \tag{2.20}$$

for all $n$. (Here we also use the obvious fact that the median of a random variable is less than twice the mean.) Choose a convergent subsequence $\{M_{n_i}\}_{i=1}^\infty$ of $\{M_n\}_{n=1}^\infty$ and set

$$\lim_{i \to \infty} M_{n_i} = \overline{M}. \tag{2.21}$$

It then follows from (2.19) and (2.21) that

$$\lim_{i \to \infty} \|G\|_{h_{n_i},p} = \overline{M} \quad \text{a.s.} \tag{2.22}$$

It follows from (2.6) and (2.17) that $\widehat{\sigma}_n^2$ is uniformly bounded. Therefore, by (2.8), for all $r > 0$,

$$E|\|G\|_{h,p} - \operatorname{med}(\|G\|_{h,p})|^r \leq C'(r), \tag{2.23}$$

for some function $C'(r)$ that depends only on $r$. We show in (2.20) that $\operatorname{med}(\|G\|_{h,p})$ is bounded uniformly in $h$. Therefore, for all $r > 0$, there exist finite constants $C(r)$ such that

$$E\|G\|_{h_n,p}^r \leq C(r) \quad \forall n \geq 1. \tag{2.24}$$

Thus, in particular, $\{\|G\|_{h_n,p}^p; \ n=1,\ldots\}$ is uniformly integrable for all $1 \leq p < \infty$. This, together with (2.22), shows that

$$\lim_{i \to \infty} E\|G\|_{h_{n_i},p}^p = \overline{M}^p. \tag{2.25}$$

Since $E\|G\|_{h_n,p}^p = (b-a)E|\eta|^p$, we have that

$$\overline{M}^p = (b-a)E|\eta|^p. \tag{2.26}$$

Thus, the bounded set $\{M_n\}_{n=1}^\infty$ has a unique limit point $\overline{M}$. It now follows from (2.19) that

$$\lim_{n \to \infty} \|G\|_{h_n,p}^p = (b-a)E|\eta|^p. \tag{2.27}$$

This gives us (2.4) when $b - a < \delta/2$. To extend the result so that it holds for any $a < b$, simply divide the interval $[a,b]$ into a finite number of subintervals with lengths $\delta/2$ and write the integral in (2.34) as a sum of integrals over these subintervals.

We now have (2.4) for fixed $a$ and $b$. Clearly, it extends to all $a$ and $b$ in a countable dense subset of $R^1$. It extends further, to all $a$ and $b$, by using the property that both the left-hand side and right-hand side of (2.27) are increasing as $a \downarrow$ and $b \uparrow$.

We conclude the proof by obtaining (2.18). Note that $\rho_h(x,y)$ is actually a function of $|x-y|$. We write $\rho_h(x,y) = \rho_h(x-y)$. Using the fact that $\rho_h$ is symmetric and setting $c = b - a$, we see that

$$\text{(2.28)} \quad \int_a^b \int_a^b |\rho_h(x-y)|\,dx\,dy = \int_0^c \int_0^c |\rho_h(x-y)|\,dx\,dy$$

$$= 2\int_0^c |\rho_h(s)|(c-s)\,ds$$

$$\text{(2.29)} \quad \leq 2(b-a) \int_0^c |\rho_h(s)|\,ds.$$

Furthermore, using the fact that $\sigma^2(h)$ is concave and monotonically increasing,

$$\sigma^2(h) \int_h^c |\rho_h(s)|\,ds$$

$$= \int_h^c (\sigma^2(s) - \sigma^2(s-h) - (\sigma^2(s+h) - \sigma^2(s)))\,ds$$

$$\text{(2.30)} \quad = \int_h^c (\sigma^2(s) - \sigma^2(s-h))\,ds - \int_{2h}^{c+h} (\sigma^2(s) - \sigma^2(s-h))\,ds$$

$$\leq \int_h^{2h} (\sigma^2(s) - \sigma^2(s-h))\,ds \leq h\sigma^2(h)$$

and

$$\sigma^2(h) \int_0^h |\rho_h(s)|\,ds$$

$$\text{(2.31)} \quad \leq \int_0^h ((\sigma^2(s+h) - \sigma^2(s)) + |\sigma^2(h-s) - \sigma^2(s)|)\,ds$$

$$\leq 2h\sigma^2(h).$$

Combining (2.28)–(2.31), we get

$$\text{(2.32)} \quad \int_a^b \int_a^b |\rho_h(x-y)|\,dx\,dy \leq 6(b-a)h,$$

which gives us (2.18). $\square$

When $G$ in Theorem 2.1 is continuous and $\sigma$ satisfies a very mild regularity condition we can take the limit in (2.4), with $h_n$ replaced by $h$.

THEOREM 2.2. *Let $G$ be as in Theorem 2.1 and assume, furthermore, that $G$ is continuous. Let $1 \leq p < \infty$ and set $h_n = 1/(\log n)^q$, where $q > p$. If*

$$\lim_{n \to \infty} \frac{\sigma(h_n - h_{n+1})}{\sigma(h_{n+1})} = 0, \tag{2.33}$$

*then*

$$\lim_{h \to 0} \int_a^b \left| \frac{G(x+h) - G(x)}{\sigma(h)} \right|^p dx = E|\eta|^p (b-a) \tag{2.34}$$

*for all $a, b \in R^1$, almost surely.*

PROOF. Without loss of generality, we assume that $b > 0$. Let

$$\|\Delta^h G\|_{p,[a,b]} := \left( \int_a^b |G(x+h) - G(x)|^p \, dx \right)^{1/p} \tag{2.35}$$

and set

$$J_G(h; a, b, p) = \frac{\|\Delta^h G\|_{p,[a,b]}}{\sigma(h)}. \tag{2.36}$$

In this notation we can write (2.4) as

$$\lim_{n \to \infty} J_G(h_n; a, b, p) = (E|\eta|^p)^{1/p} (b-a)^{1/p} \qquad \text{a.s.} \tag{2.37}$$

Fix $\delta > 0$ and consider a path for which both (2.37) holds and also the analogous statement with $b$ replaced by $2b$. We show that for such a path there exists an integer $n_1$, depending on the path and $\delta$, such that

$$|J_G(h; a, b, p) - (E|\eta|^p)^{1/p} (b-a)^{1/p}| \leq \delta \qquad \forall h \leq h_{n_1}. \tag{2.38}$$

Since we can do this for all $\delta > 0$ and all paths in a set of measure one, we get (2.34).

Set $C_0 = 2(E|\eta|^p)^{1/p}(b-a)^{1/p} \vee 1$ and $\epsilon = \delta/6C_0$. By taking $\delta$ small enough, we can assume that $\epsilon < 1/10$. Choose $N_1 > 10$ sufficiently large so that

$$\frac{\sigma(h_n - h_{n+1})}{\sigma(h_{n+1})} \leq \epsilon, \tag{2.39}$$

$$|J_G(h_n; a, b, p) - (E|\eta|^p)^{1/p}(b-a)^{1/p}| \leq \epsilon, \tag{2.40}$$

$$J_G(h_n; a, 2b, p) \leq C_0 \tag{2.41}$$

for all $n \geq N_1$. The inequality in (2.41) implies that

$$\sup_{a \leq c \leq d \leq 2b} J_G(h_n; c, d, p) \leq C_0 \qquad \forall n \geq N_1. \tag{2.42}$$

Note that for any $\zeta < h_{N_1}$ we can find an integer $m \geq N_1$ such that

(2.43) $$\zeta/2 \leq h_m \leq \zeta.$$

To see this, simply take $m = [\exp(\zeta^{-1/q})] + 1$.

To obtain (2.38), it suffices to show that it holds for all $h \in (h_{n_1+1}, h_{n_1}]$ for any $n_1 \geq N_1$. We proceed to do this. Fix $n_1$. We inductively define an increasing subsequence $\{n_j\}$, with $\lim_{j \to \infty} n_j = \infty$ beginning with $n_1$. Assume that $n_1, \ldots, n_{j-1}$, $j \geq 2$, have been defined and set $u_{j-1} := \sum_{i=1}^{j-1} h_{n_i+1}$. We take $n_j$ to be the smallest integer with

(2.44) $$h_{n_j+1} \leq h - u_{j-1}.$$

It follows from (2.43) that

(2.45) $$(h - u_{j-1})/2 \leq h_{n_j+1} \leq h - u_{j-1} < h_{n_j},$$

which implies that

(2.46) $$\lim_{j \to \infty} u_j = h.$$

It follows from the last inequality in (2.45) that $h - u_j \leq h_{n_j} - h_{n_j+1}$. Therefore, replacing $j$ by $j-1$, we have

(2.47) $$h - u_{j-1} \leq h_{n_{j-1}} - h_{n_{j-1}+1},$$

which implies, by (2.45), that

(2.48) $$h_{n_j+1} \leq h_{n_{j-1}} - h_{n_{j-1}+1}.$$

We now show that, for all $j \geq 2$,

(2.49) $$\frac{\sigma(u_j - u_{j-1})}{\sigma(u_{j-1})} \leq \epsilon^{j-1} \quad \text{and} \quad \frac{\sigma(h - u_{j-1})}{\sigma(u_{j-1})} \leq \epsilon^{j-1}.$$

To see this, we note that by (2.48) and the fact that $\sigma$ is increasing

(2.50) $$\begin{aligned}
\frac{\sigma(u_j - u_{j-1})}{\sigma(u_{j-1})} &= \frac{\sigma(h_{n_j+1})}{\sigma(u_{j-1})} \\
&= \frac{\sigma(h_{n_j+1})}{\sigma(h_{n_{j-1}+1})} \frac{\sigma(h_{n_{j-1}+1})}{\sigma(h_{n_{j-2}+1})} \cdots \frac{\sigma(h_{n_2+1})}{\sigma(u_{j-1})} \\
&\leq \frac{\sigma(h_{n_j+1})}{\sigma(h_{n_{j-1}+1})} \frac{\sigma(h_{n_{j-1}+1})}{\sigma(h_{n_{j-2}+1})} \cdots \frac{\sigma(h_{n_2+1})}{\sigma(h_{n_1+1})} \\
&\leq \frac{\sigma(h_{n_{j-1}} - h_{n_{j-1}+1})}{\sigma(h_{n_{j-1}+1})} \frac{\sigma(h_{n_{j-2}} - h_{n_{j-2}+1})}{\sigma(h_{n_{j-2}+1})} \cdots \frac{\sigma(h_{n_1} - h_{n_1+1})}{\sigma(h_{n_1+1})}.
\end{aligned}$$

The first inequality in (2.49) now follows from (2.39); the second follows similarly using (2.47).

Since (2.40) holds for all $n \geq N_1$, we have

(2.51) $$|J_G(u_1; a, b, p) - (E|\eta|^p)^{1/p}(b-a)^{1/p}| \leq \epsilon.$$

[For notational convenience, let $J_G(u_0; a, b, p) := (E|\eta|^p)^{1/p}(b-a)^{1/p}.$] For any $j \geq 1$, we have

(2.52) $$\begin{aligned}|J_G(h; a, b, p) - (E|\eta|^p)^{1/p}(b-a)^{1/p}| \\ \leq |J_G(h; a, b, p) - J_G(u_j; a, b, p)| \\ + \sum_{i=1}^{j} |J_G(u_i; a, b, p) - J_G(u_{i-1}; a, b, p)|.\end{aligned}$$

To estimate this, note that, since $\sigma$ is monotonically increasing, for any $0 < r < s$,

(2.53) $$\begin{aligned}|J_G(s; a, b, p) - J_G(r; a, b, p)| \\ = \left|\frac{\|\Delta^s G\|_{p,[a,b]}}{\sigma(s)} - \frac{\|\Delta^r G\|_{p,[a,b]}}{\sigma(r)}\right| \\ \leq \left|\frac{1}{\sigma(s)} - \frac{1}{\sigma(r)}\right| \|\Delta^r G\|_{p,[a,b]} \\ + \frac{1}{\sigma(s)} |\|\Delta^s G\|_{p,[a,b]} - \|\Delta^r G\|_{p,[a,b]}| \\ \leq \frac{|\sigma(s) - \sigma(r)|}{\sigma(r)} \frac{\|\Delta^r G\|_{p,[a,b]}}{\sigma(r)} \\ + \frac{1}{\sigma(r)} \|\Delta^s G - \Delta^r G\|_{p,[a,b]}.\end{aligned}$$

It is easy to see that the concavity of $\sigma^2$ implies the concavity of $\sigma$. Therefore, we have

(2.54) $$\frac{|\sigma(s) - \sigma(r)|}{\sigma(r)} \frac{\|\Delta^r G\|_{p,[a,b]}}{\sigma(r)} \leq \frac{\sigma(s-r)}{\sigma(r)} J_G(r; a, b, p).$$

Furthermore,

(2.55) $$\|\Delta^s G - \Delta^r G\|_{p,[a,b]} = \|\Delta^{s-r} G\|_{p,[a+r,b+r]}.$$

Consequently, for $0 < r < s$,

(2.56) $$\begin{aligned}|J_G(s; a, b, p) - J_G(r; a, b, p)| \\ \leq \frac{\sigma(s-r)}{\sigma(r)} J_G(r; a, b, p) + \frac{1}{\sigma(r)} \|\Delta^{s-r} G\|_{p,[a+r,b+r]} \\ \leq \frac{\sigma(s-r)}{\sigma(r)} (J_G(r; a, b, p) + J_G(s-r; a+r, b+r, p)).\end{aligned}$$

In particular, for any $i \geq 2$, by (2.49), we have that

$$
\begin{aligned}
&|J_G(u_i; a, b, p) - J_G(u_{i-1}; a, b, p)| \\
&\quad \leq \epsilon^{i-1}(J_G(u_{i-1}; a, b, p) + J_G(h_{n_i+1}; a+u_{i-1}, b+u_{i-1}, p)) \\
&\quad \leq \epsilon^{i-1}(J_G(u_{i-1}; a, b, p) + C_0),
\end{aligned} \tag{2.57}
$$

where, for the last step, we use (2.42).

We claim that for any $i \geq 1$

$$
J_G(u_i; a, b, p) \leq 2C_0. \tag{2.58}
$$

By (2.42), this is true for $i = 1$, without the factor of 2. However, for $i > 1$, $u_i$ need not be a member of the sequence $\{h_n\}$. To obtain (2.58), assume that it is true for all $k < i$. Then by (2.57),

$$
J_G(u_i; a, b, p) \leq C_0 + \sum_{k=2}^{i} \epsilon^{k-1} 3C_0 \leq 2C_0. \tag{2.59}
$$

It follows from (2.57) and (2.58) that

$$
|J_G(u_i; a, b, p) - J_G(u_{i-1}; a, b, p)| \leq 3\epsilon^{i-1} C_0. \tag{2.60}
$$

Using this together with (2.51) and (2.52), we see that, for any $j \geq 1$,

$$
\begin{aligned}
&|J_G(h; a, b, p) - (E|\eta|^p)^{1/p}(b-a)^{1/p}| \\
&\quad \leq |J_G(h; a, b, p) - J_G(u_j; a, b, p)| + 4\epsilon C_0.
\end{aligned} \tag{2.61}
$$

By (2.46) and the continuity of $\sigma$, we can assume that, for $j$ sufficiently large, $\sigma(u_j) \geq \sigma(h)/2$. Then using the first two lines of (2.56), (2.49) and (2.58), we see that, for all $j \geq 2$,

$$
\begin{aligned}
&|J_G(h; a, b, p) - J_G(u_j; a, b, p)| \\
&\quad \leq \frac{\sigma(h - u_j)}{\sigma(u_j)} J_G(u_j; a, b, p) \\
&\qquad + \frac{1}{\sigma(u_j)} \|\Delta^{h-u_j} G\|_{p,[a+u_j, b+u_j]} \\
&\quad \leq 2\epsilon^{j-1} C_0 + \frac{1}{\sigma(h)} \|\Delta^{h-u_j} G\|_{p,[a,2b]}.
\end{aligned} \tag{2.62}
$$

We can choose $j$ so that $h - u_j$ is arbitrarily small. Therefore, since $G$ is continuous, for a fixed path $\omega$, we can make $\|\Delta^{h-u_j} G\|_{p,[a,2b]}$ arbitrarily small. Since $\delta = 6\epsilon C_0$, we obtain (2.38). $\square$

Condition (2.33) is very weak. It is satisfied by any reasonable function one can think of, but we cannot show that it is always satisfied. In the next lemma we show that it holds when $\sigma^2(h) \geq Ch^{1/q}$, for some $q > p$. In particular, when $p = 1$, it holds for $\sigma^2(h) \geq Ch^{1-\epsilon}$ for any $\epsilon > 0$. [Since $\sigma^2$ is concave, we must have $\sigma^2(h) \geq Ch$, for some constant $C$.]

LEMMA 2.2. *When $\sigma^2(h) \geq Ch^{1/q}$, for some $q > p$, (2.33) holds.*

PROOF. Since $h_n = 1/(\log n)^q$, when $\sigma^2(h) \geq Ch^{1/q}$,

$$\sigma^2(h_n) \geq C/(\log n). \tag{2.63}$$

Suppose (2.33) does not hold. Then there exists a $\delta > 0$ and a decreasing subsequence $\{h_{n_k}\}$ of $\{h_n\}$ for which

$$\sigma(h_{n_k} - h_{n_k+1}) \geq \delta \sigma(h_{n_k+1}) \tag{2.64}$$

and $h_{n_k} - h_{n_k+1} \leq (h_{n_{k-1}} - h_{n_{k-1}+1})^2$. Using this last inequality, we see that

$$\int_{h_{n_k}-h_{n_k+1}}^{h_{n_{k-1}}-h_{n_{k-1}+1}} \frac{du}{u(\log(1/u))^{1/2}} \geq \frac{1}{4}(\log(1/(h_{n_k} - h_{n_k+1})))^{1/2}. \tag{2.65}$$

Using this, the monotonicity of $\sigma$, (2.63) and (2.64), we see that

$$\begin{aligned}
&\int_{h_{n_k}-h_{n_k+1}}^{h_{n_{k-1}}-h_{n_{k-1}+1}} \frac{\sigma(u)\,du}{u(\log(1/u))^{1/2}} \\
&\qquad \geq \frac{\delta}{4}\sigma(h_{n_k+1})(\log(1/(h_{n_k} - h_{n_k+1})))^{1/2} \\
&\qquad \geq \frac{\delta C^{1/2}}{4}\left(\frac{\log(1/(h_{n_k} - h_{n_k+1}))}{\log(n_k+1)}\right)^{1/2} > C^{1/2},
\end{aligned} \tag{2.66}$$

where for the last inequality we use the fact that, for all $n_k$ sufficiently large,

$$h_{n_k} - h_{n_k+1} \leq \frac{2q}{n_k(\log n_k)^{q+1}}. \tag{2.67}$$

Consequently, summing the left-hand side of (2.66) over all $k$ sufficiently large, we see that, for all $\alpha > 0$,

$$\int_0^\alpha \frac{\sigma(u)\,du}{u(\log(1/u))^{1/2}} = \infty. \tag{2.68}$$

This contradicts the fact that $G$ is continuous. See Example 6.4.5 in [2]. □

It is clear that the limit in (2.34) does not hold when $\sigma^2(h) = h^2$. This case includes Gaussian processes with differentiable paths. In this case

$$\lim_{h\to 0} I_G(h; a, b, p) = \int_a^b |G'(x)|^p\,dx, \tag{2.69}$$

which is not constant in general. For example, $G$ could be an integrated Brownian motion, in which case $G'$ would be the Brownian motion. Nevertheless, it is not necessary that $\sigma^2(h) \geq Ch$ for the limit to exist. We touch on this briefly in the next result for the fractional Brownian motion.

THEOREM 2.3. *Let $G$ be a fractional Brownian motion, that is, $\sigma^2(h) = h^r$, $0 < r < 2$ then (2.34) holds for all $a, b \in R^1$, almost surely.*

PROOF. Clearly, this is immediately a consequence of Theorem 2.2 for $0 < r \leq 1$, but when $1 < r < 2$, $\sigma^2(h)$ is convex. We consider this case. Let $\sigma^2(h) = h^r$, $1 < r < 2$. Analogous to (2.30), we now have

$$
\begin{aligned}
\sigma^2(h) &\int_h^c |\rho_h(s)| \, ds \\
&= \int_h^c ((\sigma^2(s+h) - \sigma^2(s)) - (\sigma^2(s) - \sigma^2(s-h))) \, ds \\
(2.70) \quad &= \int_h^c (\sigma^2(s+h) - \sigma^2(s)) \, ds - \int_0^{c-h} (\sigma^2(s+h) - \sigma^2(s)) \, ds \\
&\leq \int_{c-h}^c (\sigma^2(s+h) - \sigma^2(s)) \, ds \\
&\leq 2rc^{r-1}h^2 = 2rc^{r-1}h^{2-r}\sigma^2(h)
\end{aligned}
$$

for all $h$ sufficiently small. Also,

$$
\begin{aligned}
\sigma^2(h) &\int_0^h |\rho_h(s)| \, ds \\
(2.71) \quad &= \int_0^h ((\sigma^2(s+h) - \sigma^2(s)) + (\sigma^2(h-s) - \sigma^2(s))) \, ds \\
&\leq 2h\sigma^2(2h) \leq 8h\sigma^2(h).
\end{aligned}
$$

Consequently, when $\sigma^2(h) = h^r$, $1 < r < 2$,

$$
(2.72) \qquad \int_a^b \int_a^b |\rho_h(x-y)| \, dx \, dy \leq Ch^{2-r}.
$$

Because of the difference between (2.72) and (2.30), we must take $h_n = o(\frac{1}{(\log n)^{p/(2-r)}})$ in Lemma 2.1. This does not cause us a problem. The proof of Theorem 2.2 also works when $\sigma^2(h) = h^r$ because $\sigma$ is concave and in the proof of Theorem 2.2 the power of the $|\log h_n|$ is arbitrary. □

REMARK 2.1. Theorem 2.1, which is critical in our approach, depends on the deep Borell, Sudakov–Tsirelson Theorem. We have found a much simpler proof, based on work of Wschebor [3] that gives (2.4) for $h_n = n^{-q}$ for any $q > 2$, independent of $p$. Thus, (2.33) holds when $\sigma$ is a power. However, a sufficient condition for a Gaussian process to be continuous, when $\sigma$ is increasing, is that the integral in (2.68) is finite. This is the case, for example, if $\sigma(h) = (\log 1/h)^{-r}$ for $h \in (0, h_0]$ for some $h_0 > 0$, and $r > 1/2$. In this case (2.33) holds when $h_n = (\log n)^{-q}$, but not when $h_n = n^{-q}$.

## 3. $L^p$ moduli of continuity of squares of Gaussian processes.

The results of Section 2 immediately extend to the squares of the Gaussian processes. This is what we use to obtain results for local times.

LEMMA 3.1. *Let $\{G(x), x \in R\}$ be a mean zero continuous Gaussian process with stationary increments. Let $\sigma^2(h)$ be as defined in (2.1) and assume that*

$$(3.1) \qquad \lim_{h \to 0} \int_a^b \left| \frac{G(x+h) - G(x)}{\sigma(h)} \right|^p dx = E|\eta|^p (b - a)$$

*for all $a, b \in R^1$ almost surely, where $\eta$ is a normal random variable with mean 0 and variance 1. Then*

$$(3.2) \qquad \lim_{h \to 0} \int_a^b \left| \frac{G^2(x+h) - G^2(x)}{\sigma(h)} \right|^p dx = E|\eta|^p 2^p \int_a^b |G(x)|^p dx$$

*for all $a, b \in R^1$, almost surely.*

PROOF. Let $a = r_0 < r_1 < \cdots < r_m = b$. We have

$$
\begin{aligned}
&\int_a^b \left| \frac{G^2(x+h) - G^2(x)}{\sigma(h)} \right|^p dx \\
(3.3) \quad &= \sum_{j=1}^m \int_{r_{j-1}}^{r_j} \left| \frac{G^2(x+h) - G^2(x)}{\sigma(h)} \right|^p dx \\
&\leq 2^p \sum_{j=1}^m \int_{r_{j-1}}^{r_j} \left| \frac{G(x+h) - G(x)}{\sigma(h)} \right|^p dx \sup_{r_{j-1} \leq x \leq r_j + h} |G(x)|^p.
\end{aligned}
$$

Using (3.1), we can take the limit, as $h$ goes to zero, of the last line in (3.3) to obtain

$$
\begin{aligned}
(3.4) \quad &\limsup_{h \to 0} \int_a^b \left| \frac{G^2(x+h) - G^2(x)}{\sigma(h)} \right|^p dx \\
&\leq E|\eta|^p 2^p \sum_{j=1}^m \sup_{r_{j-1} \leq x \leq r_j} |G(x)|^p (r_j - r_{j-1}) \qquad \text{a.s.}
\end{aligned}
$$

Since $G$ has continuous sample paths, almost surely, we can take the limit of the right-hand side of (3.4), as $m$ goes to infinity and $\sup_{1 \leq j \leq m-1} r_{j+1} - r_j$ goes to zero, and use the definition of Riemann integration to get the upper bound in (3.2).

Similarly to the way we obtain (3.4), we get

$$
\begin{aligned}
(3.5) \quad & \liminf_{h \to 0} \int_a^b \left| \frac{G^2(x+h) - G^2(x)}{\sigma(h)} \right|^p dx \\
& \geq E|\eta|^p 2^p \sum_{j=1}^m \inf_{r_{j-1} \leq x \leq r_j} |G(x)|^p (r_j - r_{j-1}) \quad \text{a.s.}
\end{aligned}
$$

Taking the limit as $m$ goes to infinity and $\sup_{1 \leq j \leq m-1} r_{j+1} - r_j$ goes to zero, as in the previous paragraph, we get the lower bound in (3.2).

We have now obtained (3.2) for a fixed $a$ and $b$. We extend it to all $a, b \in R^1$ as in the proof of Theorem 2.1. $\square$

**4. Almost sure $L^p$ moduli of continuity of local times of Lévy processes.** We give some additional properties of symmetric Lévy processes $X = \{X(t), t \in R_+\}$ introduced in (1.6)–(1.10). For $0 < \alpha < \infty$ let $u^\alpha(x, y)$ denote the $\alpha$-potential density of $X$. Then

$$
(4.1) \quad u^\alpha(x, y) = \frac{1}{\pi} \int_0^\infty \frac{\cos \lambda(x-y)}{\alpha + \psi(\lambda)} d\lambda.
$$

Also, since $u^\alpha(x, y)$ is a function of $x - y$ we often write it as $u^\alpha(x - y)$.

Because of (1.9), $X$ has continuous transition probability densities, $p_t(x, y) = p_t(x - y)$; see, for example, [2], (4.74). Consequently, it is easy to see that $u^\alpha(x, y)$ is a positive definite function [2], Lemma 3.3.3. For $0 < \alpha < \infty$, let

$$
\begin{aligned}
(4.2) \quad \sigma_\alpha^2(x - y) & := u^\alpha(x, x) + u^\alpha(y, y) - 2u^\alpha(x, y) \\
& = 2(u^\alpha(0) - u^\alpha(x - y)) \\
& = \frac{4}{\pi} \int_0^\infty \sin^2 \frac{\lambda(x-y)}{2} \frac{1}{\alpha + \psi(\lambda)} d\lambda.
\end{aligned}
$$

We can also consider $u^\alpha(x, y)$, $0 < \alpha < \infty$, as the covariance of a mean zero stationary Gaussian process, which we denote by $G_\alpha = \{G_\alpha(x), x \in R\}$. We have

$$
(4.3) \quad E(G_\alpha(x) - G_\alpha(y))^2 = \sigma_\alpha^2(x - y).
$$

Note that the covariance of $G_\alpha$ is the 0-potential density of a Lévy process killed at the end of an independent exponential time with mean $1/\alpha$. Thus, $G_\alpha$ is an associated Gaussian process in the nomenclature of [2].

We are interested in those Lévy processes with 1-potential density given by (4.1) for which the stationary Gaussian processes $G_1$, defined by (4.3), are continuous and satisfy (3.1). We refer to these processes as Lévy processes of class A. Since the Gaussian processes $G_1$ are continuous, we know that the Lévy processes of class A have jointly continuous local times ([2], Theorem 9.4.1, (1)).

We now use the Eisenbaum Isomorphism Theorem, as employed in [2], Theorem 10.4.1, to obtain the following $L^p$ moduli of continuity for the local times of these Lévy processes.

LEMMA 4.1. *Let $X = \{X(t), t \in R_+\}$ be a real-valued symmetric Lévy process of class A with 1-potential density $u^1(x, y)$ and let $\{L_t^x, (t, x) \in R_+ \times R\}$ be the local time of $X$. Then, for almost all $t \in R_+$,*

$$(4.4) \qquad \lim_{h \downarrow 0} \int_a^b \left| \frac{L_t^{x+h} - L_t^x}{\sigma_1(h)} \right|^p dx = 2^{p/2} E|\eta|^p \int_a^b |L_t^x|^{p/2} dx$$

*for all $a, b \in R^1$, almost surely.*

PROOF. By Lemma 3.1,

$$(4.5) \qquad \begin{aligned} \lim_{h \to 0} \int_a^b & \left| \frac{G_1^2(x+h)/2 - G_1^2(x)/2}{\sigma_1(h)} \right|^p dx \\ &= 2^{p/2} E|\eta|^p \int_a^b |G_1^2(x)/2|^{p/2} dx \end{aligned}$$

for all $a, b \in R^1$ almost surely, where $\eta$ is a normal random variable with mean 0 and variance 1. A simple modification of the proof of Lemma 3.1 shows that, for all $s$,

$$(4.6) \qquad \begin{aligned} \lim_{h \to 0} \int_a^b & \left| \frac{(G_1(x+h) + s)^2/2 - (G_1(x) + s)^2/2}{\sigma_1(h)} \right|^p dx \\ &= 2^{p/2} E|\eta|^p \int_a^b |(G_1(x) + s)^2/2|^{p/2} dx \end{aligned}$$

for all $a, b \in R^1$ almost surely.

Let $\omega \in \Omega_{G_1}$ denote the probability space of $G_1$ and fix $\omega \in \Omega_{G_1}$. Using the notation of (2.7),

$$(4.7) \qquad \begin{aligned} & \|L_t + (G_1(\omega) + s)^2/2\|_{h,p}^p \\ &= \int_a^b \left| \frac{(L_t^{x+h} - L_t^x + (G_1(x+h, \omega) + s)^2/2 - (G_1(x, \omega) + s)^2/2)}{\sigma_1(h)} \right|^p dx. \end{aligned}$$

It follows from the Eisenbaum Isomorphism Theorem that, for any $s \neq 0$, an almost sure event for $(G_1(\omega) + s)^2/2$ is also an almost sure event for $L_t^{\cdot} + (G_1(\omega) + s)^2/2$, for almost all $t \in R_+$; see [2], Lemma 9.1.2. (Here $X$ and $G_1$ are independent.) Therefore, (4.6) implies that, for almost all $\omega \in \Omega_{G_1}$ and for almost all $t \in R_+$,

$$\lim_{h \downarrow 0} \|L_t + (G_1(\omega) + s)^2/2\|_{h,p}$$

(4.8)
$$= 2^{1/2}(E|\eta|^p)^{1/p} \left( \int_a^b |L_t^x + (G_1(x,\omega) + s)^2/2|^{p/2} \, dx \right)^{1/p}$$

for all $a, b \in R^1$ almost surely (with respect to $\Omega_X$). Consequently, for almost all $\omega \in \Omega_{G_1}$ and for almost all $t \in R_+$,

$$\limsup_{h \downarrow 0} \|L_t\|_{h,p}$$

(4.9)
$$\leq 2^{1/2}(E|\eta|^p)^{1/p}$$
$$\times \left( \left( \int_a^b |L_t^x|^{p/2} \, dx \right)^{1/p} + \left( \int_a^b |(G_1(x,\omega) + s)^2/2|^{p/2} \, dx \right)^{1/p} \right)$$
$$+ \limsup_{h \downarrow 0} \int_a^b \left| \frac{(G_1(x+h,\omega) + s)^2/2 - (G_1(x,\omega) + s)^2/2}{\sigma_1(h)} \right|^p dx$$

for all $a, b \in R^!$ almost surely. Using (4.6) on the last term in (4.9), we see that, for almost all $\omega \in \Omega_{G_1}$ and for almost all $t \in R_+$,

$$\limsup_{h \downarrow 0} \|L_t\|_{h,p} \leq 2^{1/2}(E|\eta|^p)^{1/p}$$

(4.10)
$$\times \left( \left( \int_a^b |L_t^x|^{p/2} \, dx \right)^{1/p} + 2 \left( \int_a^b |(G_1(x,\omega) + s)^2/2|^{p/2} \, dx \right)^{1/p} \right)$$

for all $a, b \in R^1$ almost surely. And since this holds for all $s \neq 0$, we get that, for almost all $\omega \in \Omega_{G_1}$ and for almost all $t \in R_+$,

$$\limsup_{h \downarrow 0} \|L_t\|_{h,p}$$

(4.11)
$$\leq 2^{1/2}(E|\eta|^p)^{1/p}$$
$$\times \left( \left( \int_a^b |L_t^x|^{p/2} \, dx \right)^{1/p} + 2 \left( \int_a^b |G_1^2(x,\omega)/2|^{p/2} \, dx \right)^{1/p} \right)$$

for all $a, b \in R^1$ almost surely.

Since $G_1$ has continuous sample paths, it follows from [2], Lemma 5.3.5, that, for all $\epsilon > 0$,

(4.12)
$$P\left( \sup_{x \in [a,b]} |G_1(x)| \leq \epsilon \right) > 0.$$

Therefore, we can choose $\omega$ in (4.11) so that the integral involving the Gaussian process can be made arbitrarily small. Thus, for almost all $t \in R^1$,

$$(4.13) \qquad \limsup_{h \downarrow 0} \|L_t\|_{h,p} \leq 2^{1/2} (E|\eta|^p)^{1/p} \left( \int_a^b |L_t^x|^{p/2} \, dx \right)^{1/p}$$

for all $a, b \in R^1$, almost surely. By the same methods, we can obtain the reverse of (4.13) for the limit inferior. $\square$

Analogous to the definition of $\sigma_\alpha^2$ in (4.2), we set

$$(4.14) \qquad \begin{aligned} \sigma_0^2(x) &:= \lim_{\alpha \to 0} 2(u^\alpha(0) - u^\alpha(x)) \\ &= \frac{4}{\pi} \int_0^\infty \sin^2 \frac{\lambda x}{2} \frac{1}{\psi(\lambda)} \, d\lambda. \end{aligned}$$

By (1.9) and the fact that $\lambda^2 = O(\psi(\lambda))$ as $\lambda \to 0$ (see [2], (4.72) and (4.77)), the integral in (4.14) is finite, so that $\sigma_0$ is well defined whether or not $X$ has a 0-potential density.

For later reference, we note that by the definition of the $\alpha$-potential density of $X$ and (4.14)

$$(4.15) \qquad \begin{aligned} \sigma_0^2(x) &= 2 \lim_{\alpha \to 0} \int_0^\infty e^{-\alpha t} (p_t(0) - p_t(x)) \, dt \\ &= 2 \int_0^\infty (p_t(0) - p_t(x)) \, dt. \end{aligned}$$

Lemma 4.1 is very close Theorem 1.2. However, Lemma 4.1 requires that $G_1$ satisfies (3.1). Theorem 2.2, which gives conditions for Gaussian processes to satisfy (3.1), requires that $\sigma_1^2$ is concave at the origin. It is easier to verify concavity for $\sigma_0^2$. That is why we use $\sigma_0^2$ in Theorem 1.2. We proceed to use Lemma 4.1 and some observations about $\sigma_1^2$ and $\sigma_0^2$ to prove Theorem 1.2.

We need some general facts about Gaussian processes with stationary increments. Let $\mu$ be a measure on $(0, \infty)$ that satisfies (1.8). Let

$$(4.16) \qquad \phi(x) := \frac{4}{\pi} \int_0^\infty \sin^2 \frac{\lambda x}{2} \, d\mu(\lambda).$$

The function $\phi(x)$ determines a mean zero Gaussian process with stationary increments $H = \{H(x), x \in R^1\}$ with $H(0) = 0$, by the relationship

$$(4.17) \qquad E((H(x) - H(y))^2) = \phi(x - y).$$

(This is because it follows from (4.17) that

$$(4.18) \qquad EH(x)H(y) = \tfrac{1}{2}(\phi(x) + \phi(y) - \phi(x - y)).$$

It is easy to see that $EH(x)H(y)$ is positive definite and, hence, determines a mean zero Gaussian process; see, e.g., [2], 5.252.)

We consider three such Gaussian processes, $G_0$, and $\overline{G}_\alpha$ and $\widetilde{G}_\alpha$ for $\alpha > 0$, determined by

$$\sigma_0^2(h) = \frac{4}{\pi} \int_0^\infty \sin^2 \frac{\lambda h}{2} \frac{1}{\psi(\lambda)} \, d\lambda, \tag{4.19}$$

$$\overline{\sigma}_\alpha^2(h) = \frac{4}{\pi} \int_0^\infty \sin^2 \frac{\lambda h}{2} \frac{1}{(\alpha + \psi(\lambda))} \, d\lambda, \tag{4.20}$$

$$\widetilde{\sigma}_\alpha^2(h) = \frac{4}{\pi} \int_0^\infty \sin^2 \frac{\lambda h}{2} \frac{\alpha}{\psi(\lambda)(\alpha + \psi(\lambda))} \, d\lambda, \tag{4.21}$$

as described in the previous paragraph. Note that $\overline{G}_\alpha(x) = G_\alpha(x) - G_\alpha(0)$, $x \in R^1$, for $G_\alpha$ as defined in (4.3). Therefore, the increments of $\overline{G}_\alpha$ and $G_\alpha$ are the same and, $\overline{\sigma}_\alpha^2 = \sigma_\alpha^2$, defined in (4.3).

Obviously,

$$\sigma_0^2(h) = \overline{\sigma}_\alpha^2(h) + \widetilde{\sigma}_\alpha^2(h). \tag{4.22}$$

Let $\overline{G}_\alpha$ and $\widetilde{G}_\alpha$ be independent. It follows from (4.22) that $\overline{G}_\alpha + \widetilde{G}_\alpha$ is a version of $G_0$. In this sense we can write

$$G_0(x) = \overline{G}_\alpha(x) + \widetilde{G}_\alpha(x), \qquad x \in R^1. \tag{4.23}$$

We show in [2], Lemma 7.4.8, that

$$\lim_{h \to 0} \frac{\sigma_0(h)}{\sigma_\alpha(h)} = 1. \tag{4.24}$$

This shows that $G_0$ has continuous paths if and only if $\overline{G}_\alpha$, or equivalently, $G_\alpha$, has continuous paths. Furthermore, by (4.22) and (4.24), if $G_\alpha$ has continuous paths, so does $\widetilde{G}_\alpha$. (These facts about continuity follows from [2], Lemma 5.5.2 and Theorem 5.3.10. See also [1], Chapter 15, Section 3.)

LEMMA 4.2. *Let $\sigma_0$, $\widetilde{\sigma}_\alpha$ and $\psi(\lambda)$ be as given in (4.19) and (4.21) and assume that $\psi(\lambda)$ satisfies (1.12). Assume also that $h^{2-\gamma'} = O(\sigma_0^2(h))$ for some $\gamma' > 0$ as $h \downarrow 0$. Then for all $\alpha > 0$, there exists an $\epsilon > 0$ such that*

$$\widetilde{\sigma}_\alpha^2(h) = O(h^\epsilon \sigma_0^2(h)) \qquad as\ h \downarrow 0. \tag{4.25}$$

PROOF. Let $\delta = \gamma'/4 < 1$. By (1.12), there exists an $M \in R^1$ such that $\psi(\lambda) \geq \lambda^\gamma$ for all $\lambda \geq M \vee 1$. Then

$$\widetilde{\sigma}_\alpha^2(h) \leq \frac{h^2}{\pi} \left( \int_0^M \frac{\lambda^2}{\psi(\lambda)} \, d\lambda + \int_M^{(1/h)^\delta} \lambda^2 \, d\lambda \right)$$

$$+ \frac{\alpha}{\inf_{x \geq (1/h)^\delta}(\alpha + \psi(x))} \int_{(1/h)^\delta}^\infty \sin^2 \frac{\lambda h}{2} \frac{1}{\psi(\lambda)} \, d\lambda \tag{4.26}$$

$$\leq 0(h^{2-3\gamma'/4}) + 0(h^{\delta\gamma} \sigma_0^2(h))$$

which implies (4.25). (Here we use the fact that $\lambda^2/\psi(\lambda)$ is bounded on $[0, M]$; see, e.g., [2], Lemma 4.2.2.) □

PROOF OF THEOREM 1.2. In this section we prove this theorem with "all $t \in R_+$" replaced by "almost all $t \in R_+$." We complete the proof of this theorem in Section 5.

Since $L$ has continuous local times, it follows from [2], Theorem 9.4.1, (1), that $G_1$, the stationary Gaussian process with covariance $u_1$, is continuous almost surely. Therefore, by the remarks made prior to the statement of Lemma 4.2, $G_1$, $\overline{G}_1$, $\widetilde{G}_1$ and $G_0$ are all continuous almost surely.

Using (4.23), we see that

$$
(4.27) \quad \left| \left( \int_a^b \left| \frac{G_0(x+h) - G_0(x)}{\sigma_0(h)} \right|^p dx \right)^{1/p} - \left( \int_a^b \left| \frac{G_1(x+h) - G_1(x)}{\sigma_0(h)} \right|^p dx \right)^{1/p} \right|
$$
$$
\leq \left( \int_a^b \left| \frac{\widetilde{G}_1(x+h) - \widetilde{G}_1(x)}{\sigma_0(h)} \right|^p dx \right)^{1/p}.
$$

We show below that the last integral in (4.27) goes to zero as $h \downarrow 0$. Furthermore, by Theorem 2.2, the limit of the first integral in (4.27) goes to $E|\eta|^p(b-a)$ almost surely as $h \downarrow 0$. Consequently, the limit of the second integral in (4.27) also goes to $E|\eta|^p(b-a)$ almost surely as $h \downarrow 0$. Using (4.24), we have

$$
(4.28) \quad \lim_{h \to 0} \int_a^b \left| \frac{G_1(x+h) - G_1(x)}{\sigma_1(h)} \right|^p dx = E|\eta|^p(b-a) \quad \text{a.s.}
$$

This shows that $X$ is a Lévy process of class A (see page 17), so (4.4) holds. Using (4.24) again gives (1.13).

Note that by (4.25) there exists an $\epsilon > 0$ such that

$$
(4.29) \quad \widetilde{\sigma}_1^2(h) \leq h^\epsilon \sigma_0^2(h) \quad \text{for } h \in [0, h_0]
$$

for some $h_0 > 0$. Therefore, by [2], Theorem 7.2.1,

$$
(4.30) \quad C(h^\epsilon \sigma_0^2(h) \log 1/h)^{1/2}
$$

is a uniform modulus of continuity for $\widetilde{G}_\alpha$. It follows from this that the last integral in (4.27) goes to zero as $h \downarrow 0$. □

REMARK 4.1. The simplest and perhaps most important application of Theorem 1.2 is to symmetric stable processes with index $1 < \beta \leq 2$. In this case $\psi(\lambda) = |\lambda|^\beta$. (Stable processes with index $\beta \leq 1$ do not have local times.) By a change of variables, we see that

$$\sigma_0^2(h) = h^{\beta-1} \frac{4}{\pi} \int_0^\infty \left(\sin^2 \frac{s}{2}\right) \frac{1}{|s|^\beta} ds$$

(4.31)

$$= h^{\beta-1} \frac{1}{\Gamma(\beta)\sin((\pi/2)(\beta-1))}.$$

The calculation that gives the last line is given in [2], (4.94) and (4.99)–(4.102), however, note that the numerator in [2], (4.102), should be one.

When $\beta = 2$ the Lévy process is $\{\sqrt{2}B_t, t \in R_+\}$, where $\{B_t, t \in R_+\}$ is a standard Brownian motion. The factor $\sqrt{2}$ occurs because the Lévy exponent in this case is $\lambda^2$ rather than $\lambda^2/2$.

PROOF OF THEOREM 1.1. This is an immediate application of Theorem 1.2 in which we calculate (1.10) with $\psi(\lambda) = \lambda^2/2$. Thus, the function $\sigma_0^2(h)$ for the Brownian motion is twice the last line in (4.31), which in this case is simply $2h$. □

We have a much larger class of concrete examples to which we can apply Theorem 1.2. In [2], Section 9.6, we consider a case of Lévy processes which we call stable mixtures. Using stable mixtures, we show in [2], Corollary 9.6.5, that for any $0 < \beta < 1$ and function $g$ which is regularly varying at infinity with positive index or is slowly varying at infinity and increasing, there exists a Lévy process for which the corresponding function $\sigma_0^2(h)$ is concave and satisfies

(4.32) $$\sigma_0^2(h) \sim |h|^\beta g(\log 1/|h|) \qquad \text{as } h \to 0.$$

Moreover, if in addition,

(4.33) $$\int_0^1 \frac{dx}{g(x)} < \infty,$$

the above statement is also valid when $\beta = 1$. Since $\sigma_0^2$ is regularly varying, (2.33) holds. Also, in [2], Section 9.6, the characteristic exponents of stable mixtures is given explicitly and it is easy to see that they satisfy (1.12).

**5. Convergence in $L^m$.** In Section 4 Theorem 1.2 is only proved for almost every $t$ (see page 21). To obtain Theorem 1.2 for all $t$, we need additional information which is contained in the next theorem. This theorem is also interesting on its own.

THEOREM 5.1. *Under the hypotheses of Theorem 1.2,*

(5.1) $$\lim_{h \downarrow 0} \int_a^b \left| \frac{L_t^{x+h} - L_t^x}{\sigma_0(h)} \right|^p dx = 2^{p/2} E|\eta|^p \int_a^b |L_t^x|^{p/2} dx$$

in $L^m$ *uniformly in t on any bounded interval of* $R_+$, *for all* $m \geq 1$.

The proof follows from several lemmas on moments of the $L^m$ norm of various functions of the local times. We begin with a formula for the moments of local times. For a proof, see [2], Lemma 10.5.5.

LEMMA 5.1. *Let* $X = \{X(t), t \in R_+\}$ *be a symmetric Lévy process and let* $\{L_t^x, (t,x) \in R_+ \times R\}$ *be the local times of* $X$. *Then for all* $x, y, z \in R$, $t \in R_+$ *and integers* $m \geq 1$,

$$(5.2) \quad E^z((L_t^x)^m) = m! \int \cdots \int_{0 < t_1 < \cdots < t_m < t} p_{t_1}(x-z) \prod_{i=2}^m p_{\Delta t_i}(0) \prod_{i=1}^m dt_i,$$

*where* $p_t$ *is the probability density function of* $X(t)$ *and* $\Delta t_i = t_i - t_{i-1}$.

*Furthermore,*

$$
\begin{aligned}
&E^z((L_t^x - L_t^y)^{2m}) \\
(5.3) \quad &= (2m)! \int \cdots \int_{0 < t_1 < \cdots < t_{2m} < t} (p_{t_1}(x-z) + p_{t_1}(y-z)) \\
&\qquad \times \prod_{i=2}^{2m}(p_{\Delta t_i}(0) - (-1)^{2m-i} p_{\Delta t_i}(x-y)) \prod_{i=1}^m dt_i.
\end{aligned}
$$

Let $Z$ be a random variable on the probability space of $X$. We denote the $L^m$ norm of $Z$ with respect to $P^0$ by $\|Z\|_m$. Let

$$(5.4) \quad V(t) = \int_0^t p_s(0) \, ds.$$

The next lemma follows easily from Lemma 5.1 and the fact that $p_s(x) \leq p_s(0)$ for all $x \in R$, and uses the representation of $\sigma_0$ in the last line of (4.15). For (5.6), we also use the fact that $L_t^x - L_s^x = L_{t-s}^x \circ \theta_s$ together with the Markov property.

LEMMA 5.2. *Let* $X = \{X(t), t \in R_+\}$ *be a real-valued symmetric Lévy process and let* $\{L_t^x, (t,x) \in R_+ \times R\}$ *be the local times of* $X$. *Then for all* $x, y \in R$, $s, t \in R_+$ *and integers* $m \geq 1$,

(5.5) $\qquad \|L_t^x - L_t^y\|_{2m} \leq C(m) V^{1/2}(t) \sigma_0(x-y),$

(5.6) $\qquad \|L_t^x - L_s^x\|_m \leq C'(m) V(t-s),$

(5.7) $\qquad \|L_t^x\|_m \leq C'(m) V(t),$

*where* $C(m)$ *and* $C'(m)$ *are constants depending only on* $\beta$ *and* $m$.

It is clear that the inequality in (5.6) is unchanged if we take the norm with respect to $P^z$, for any $z \in R$. The same observation applies to (5.5) since it only depends on $|x - y|$.

In the next lemma we use notation introduced in (2.7), except that $\sigma$ is replaced by $\sigma_0$.

LEMMA 5.3. *Let $X = \{X(t), t \in R_+\}$ be a real valued symmetric Lévy process and let $\{L_t^x, (t,x) \in R_+ \times R\}$ be the local times of $X$. Then for all $h > 0$, $s, t \in R_+$, with $s \leq t$, $p \geq 1$ and integers $m \geq 1$,*

$$(5.8) \quad \|\|L_t\|_{h,p}^p - \|L_s\|_{h,p}^p\|_m \leq C(p,m) V^{(p-1)/2}(t) V^{1/2}(t-s)(b-a).$$

*In particular,*

$$(5.9) \quad \|\|L_t\|_{h,p}^p\|_m^{1/p} \leq C'(p,m) V^{1/2}(t)(b-a)^{1/p},$$

*where $C(p,m)$ and $C'(p,m)$ are constants depending only on $p$ and $m$.*

*Similarly, for any $r \geq 1$,*

$$(5.10) \quad \left\|\int_a^b |L_t^x|^r \, dx - \int_a^b |L_s^x|^r \, dx\right\|_m \leq D(r,m) V^{r-1}(t) V(t-s)(b-a).$$

*In particular,*

$$(5.11) \quad \left\|\int_a^b |L_t^x|^r \, dx\right\|_m \leq D'(r,m) V^r(t)(b-a).$$

*For any $0 < r \leq 1$,*

$$(5.12) \quad \left\|\int_a^b |L_t^x|^r \, dx - \int_a^b |L_s^x|^r \, dx\right\|_m \leq D(r,m) V^r(t-s)(b-a).$$

*In particular,*

$$(5.13) \quad \left\|\int_a^b |L_t^x|^r \, dx\right\|_m \leq D'(r,m) V^r(t)(b-a),$$

*where $D(r,m)$ and $D'(r,m)$ are constants depending only on $r$ and $m$.*

PROOF. Set

$$(5.14) \quad \Delta^h L_t^x = L_t^{x+h} - L_t^x.$$

Suppose that $u \geq v \geq 0$. Writing $u^p - v^p$ as the integral of its derivative, we see that

$$(5.15) \quad u^p - v^p \leq p(u-v)u^{p-1}.$$

Therefore, it follows from (5.15) and the Schwarz inequality that

$$\|\|L_t\|_{h,p}^p - \|L_s\|_{h,p}^p\|_m$$

(5.16)
$$\leq \int_a^b \frac{1}{\sigma_0^p(h)} \||\Delta^h L_t^x|^p - |\Delta^h L_s^x|^p\|_m \, dx$$

$$\leq \int_a^b \frac{p}{\sigma_0^p(h)} (\||\Delta^h L_t^x|^{p-1}\|_{2m} + \||\Delta^h L_s^x|^{p-1}\|_{2m})$$

$$\times \|\Delta^h L_t^x - \Delta^h L_s^x\|_{2m} \, dx.$$

Let $r$ be the smallest even integer greater than or equal to $2m(p-1)$. Then by Hölder's inequality and (5.5), we see that

(5.17)
$$\||\Delta^h L_t^x|^{p-1}\|_{2m} \leq \|\Delta^h L_t^x\|_r^{p-1}$$

$$\leq D(m) V^{(p-1)/2}(t) \sigma_0^{p-1}(h),$$

where $D(m) = (C(r))^{p-1}$ and $C(r)$ is the constant in (5.5). (Clearly, this inequality also holds with $t$ replaced by any $s \leq t$.)

It follows from (5.5) and the remark immediately following the statement of Lemma 5.2, that for all $z \in R$,

(5.18)
$$(E^z (\Delta^h L_{t-s}^x)^{2m})^{1/2m} = \|\Delta^h L_{t-s}^{x-z}\|_{2m}$$

$$\leq C(m) V^{1/2}(t-s) \sigma_0(h).$$

Consequently,

(5.19)
$$\|\Delta^h L_t^x - \Delta^h L_s^x\|_{2m} = \|\Delta^h L_{t-s}^x \circ \theta_s\|_{2m}$$

$$= (E^0 \{E^{X_s} (\Delta^h L_{t-s}^x)^{2m}\})^{1/2m}.$$

$$\leq C(m) V^{1/2}(t-s) \sigma_0(h).$$

It follows from (5.16), (5.17) and (5.19), and the fact that $s \leq t$, that

(5.20)
$$\|\|L_t\|_{h,p}^p - \|L_s\|_{h,p}^p\|_m$$

$$\leq 2p D(m,p) C(m) V^{(p-1)/2}(t) V^{1/2}(t-s)(b-a).$$

This gives (5.8). The statement in (5.9) follows from (5.8) by setting $s = 0$.

To prove (5.10), we take $s < t$, and note that

(5.21)
$$\left\| \int_a^b |L_t^x|^r \, dx - \int_a^b |L_s^x|^r \, dx \right\|_m$$

$$\leq \int_a^b \||L_t^x|^r - |L_s^x|^r\|_m \, dx \leq (b-a) \sup_x \||L_t^x|^r - |L_s^x|^r\|_m.$$

It follows from (5.15) with $p$ replaced by $r \geq 1$, followed by the Cauchy–Schwarz inequality, that

$$(5.22) \qquad \||L_t^x|^r - |L_s^x|^r\|_m \leq r\|L_t^x - L_s^x\|_{2m} \||L_t^x|^{r-1}\|_{2m}.$$

As in (5.17), we have

$$(5.23) \qquad \||L_t^x|^{r-1}\|_{2m} \leq \|L_t^x\|_q^{r-1},$$

where $q$ is the smallest even integer greater than or equal to $2m(r-1)$. The inequality in (5.10) now follows from (5.6) and (5.7). The inequality in (5.11) follows from (5.10) by setting $s = 0$.

When $0 \leq r \leq 1$ we have

$$(5.24) \qquad 0 \leq |L_t^x|^r - |L_s^x|^r \leq |L_t^x - L_s^x|^r,$$

so that

$$(5.25) \qquad \||L_t^x|^r - |L_s^x|^r\|_m \leq \||L_t^x - L_s^x|^r\|_m \leq \|L_t^x - L_s^x\|_q^r,$$

where $q$ is the smallest integer greater than or equal to $rm$. The inequality in (5.12) now follows from (5.6). The inequality in (5.13) follows from (5.12) by setting $s = 0$. $\square$

PROOF OF THEOREM 5.1. Although it is usually easier to prove convergence in $L^m$ than it is to prove convergence almost surely, the only way that we know to prove this theorem is by using Theorem 1.2. Fix $a < b$. For $h > 0$, let

$$(5.26) \qquad H_h(t) = \int_a^b \frac{|L_t^{x+h} - L_t^x|^p}{\sigma_0^p(h)}\, dx - 2^{p/2} E|\eta|^p \int_a^b |L_t^x|^{p/2}\, dx.$$

It follows from Theorem 1.2 and Fubini's theorem that there exists dense subset $D \subseteq R^+$, such that, for each $t \in D$, $H_h(t)$ converges to 0 almost surely.

By (5.9) and (5.11), we have that, for any $m$,

$$(5.27) \qquad \|H_h(t)\|_m \leq C(m, b-a, t) < \infty,$$

where the function $C(m, b-a, t)$ is independent of $h$. In particular, for each $t$, the collection $\{H_h(t); h > 0\}$ is uniformly integrable. Consequently, for any $m \geq 1$,

$$(5.28) \qquad \lim_{h \downarrow 0} \|H_h(t)\|_m = 0 \qquad \forall t \in D.$$

Fix $T > 0$. By (5.8), (5.10) and (5.12) for any $m \geq 1$ and any $\epsilon > 0$, we can find a $\delta > 0$ such that

$$(5.29) \qquad \sup_{\substack{0 \leq s, t \leq T \\ |s-t| \leq \delta}} \|H_h(s) - H_h(t)\|_m \leq \epsilon \qquad \forall h > 0.$$

Choose a finite set $\{t_1,\ldots,t_k\}$ in $D \cap [0,T]$ such that $\bigcup_{j=1}^{k}[t_j - \delta, t_j + \delta]$ covers $[0,T]$. By (5.28), we can choose an $h_\epsilon$ such that

$$\sup_{j=1,\ldots,k} \|H_h(t_j)\|_m \leq \epsilon \quad \forall h \leq h_\epsilon. \tag{5.30}$$

Combined with (5.29), this shows that

$$\sup_{0 \leq s \leq T} \|H_h(s)\|_m \leq 2\epsilon \quad \forall h \leq h_\epsilon. \tag{5.31}$$

$\square$

Using Theorem 5.1 we can now complete the proof of Theorem 1.2.

PROOF OF THEOREM 1.2 CONTINUED. Fix $-\infty < a < b < \infty$. What we have already proved (see page 21) implies that we can find a dense subset $T' \in R_+$ such that

$$\lim_{h \downarrow 0} \int_a^b \left| \frac{L_s^{x+h} - L_s^x}{\sigma_0(h)} \right|^p dx = 2^{p/2} E|\eta|^p \int_a^b |L_s^x|^{p/2}\, dx \tag{5.32}$$

for all $s \in T'$ almost surely. Fix $t > 0$, and let $s_n, n=1,\ldots$, be a sequence in $T'$ with $s_n \uparrow t$. Using the additivity of local times, we have

$$\Delta^h L_t^x - \Delta^h L_{s_n}^x = \Delta^h L_{t-s_n}^x \circ \theta_{s_n}, \tag{5.33}$$

so that, in the notation of (2.35),

$$A_n := \limsup_{h \downarrow 0} \frac{1}{\sigma_0(h)} \left| \|\Delta^h L_t^x\|_{p,[a,b]} - \|\Delta^h L_{s_n}^x\|_{p,[a,b]} \right|$$

$$\leq \limsup_{h \downarrow 0} \frac{1}{\sigma_0(h)} \|\Delta^h L_t^x - \Delta^h L_{s_n}^x\|_{p,[a,b]} \tag{5.34}$$

$$= \limsup_{h \downarrow 0} \frac{1}{\sigma_0(h)} \|\Delta^h L_{t-s_n}^x \circ \theta_{s_n}\|_{p,[a,b]}.$$

Let $\bar{X}_r = X_{r+s_n} - X_{s_n}, r \geq 0$. Note that $\{\bar{X}_r; r \geq 0\}$ is a copy of $\{X_r; r \geq 0\}$ that is independent of $X_{s_n}$. Let $\{\bar{L}_r^x; (x,r) \in R^1 \times R_+\}$ denote the local time for the process $\{\bar{X}_r; r \geq 0\}$. It is easy to check that

$$L_{t-s_n}^x \circ \theta_{s_n} = \bar{L}_{t-s_n}^{x-X_{s_n}}. \tag{5.35}$$

Therefore,

$$\|\Delta^h L_{t-s_n}^x \circ \theta_{s_n}\|_{p,[a,b]} = \|\Delta^h \bar{L}_{t-s_n}^{x-X_{s_n}}\|_{p,[a,b]} \tag{5.36}$$

$$= \|\Delta^h \bar{L}_{t-s_n}^x\|_{p,[a-X_{s_n},b-X_{s_n}]}.$$

Since $X_{s_n}$ is independent of $\{\bar{X}_r; r \geq 0\}$, it follows from Theorem 5.1 that, conditional on $X_{s_n}$,

$$\lim_{h \downarrow 0} \frac{1}{\sigma_0(h)} \|\Delta^h \bar{L}^x_{t-s_n}\|_{p,[a-X_{s_n},b-X_{s_n}]}$$
$$(5.37)$$
$$= 2^{1/2}(E|\eta|^p)^{1/p} \|\bar{L}^x_{t-s_n}\|^{1/2}_{p/2,[a-X_{s_n},b-X_{s_n}]} \quad \text{in } L^1_{\bar{X}},$$

where $L^1_{\bar{X}}$ denotes $L^1$ with respect to $\bar{X}$.

We now use (5.37) followed by Hölder's inequality, and then either (5.11) for $1 \leq p/2 < \infty$, or (5.13) for $0 < p/2 < 1$, to see that

$$E(A_n|X_{s_n}) \leq 2^{1/2}(E|\eta|^p)^{1/p} E(\|\bar{L}^x_{t-s_n}\|^{1/2}_{p/2,[a-X_{s_n},b-X_{s_n}]} | X_{s_n})$$
$$(5.38) \qquad \leq 2^{1/2}(E|\eta|^p)^{1/p} |E(\|\bar{L}^x_{t-s_n}\|^{p/2}_{p/2,[a-X_{s_n},b-X_{s_n}]} | X_{s_n})|^{1/p}$$
$$\leq 2^{1/2}(E|\eta|^p D'(\beta,p/2,1)(b-a))^{1/p} V^{1/2}(t-s_n).$$

Therefore,

$$(5.39) \qquad E(A_n) \leq CV^{1/2}(t-s_n),$$

where $C < \infty$, is independent of $n$. Since $T'$ is dense in $R^+$, we can choose a sequence $\{s_n\} \in T'$, so that $\sum_{n=1}^{\infty} V^{1/2}(t-s_n) < \infty$. Therefore, by (5.39) and the Borel–Cantelli Lemma,

$$(5.40) \qquad \lim_{n \to \infty} A_n = 0 \quad \text{a.s.}$$

The proof of this theorem is completed by observing that, for each $n$,

$$\limsup_{h \downarrow 0} \frac{1}{\sigma_0(h)} \|\Delta^h L^x_t\|_{p,[a,b]} \leq \limsup_{h \downarrow 0} \frac{1}{\sigma_0(h)} \|\Delta^h L^x_{s_n}\|_{p,[a,b]} + A_n$$
$$= 2^{1/2}(E|\eta|^p)^{1/p} \|L^x_{s_n}\|^{1/2}_{p/2,[a,b]} + A_n,$$
$$\liminf_{h \downarrow 0} \frac{1}{\sigma_0(h)} \|\Delta^h L^x_t\|_{p,[a,b]} \geq \liminf_{h \downarrow 0} \frac{1}{\sigma_0(h)} \|\Delta^h L^x_{s_n}\|_{p,[a,b]} - A_n$$
$$= 2^{1/2}(E|\eta|^p)^{1/p} \|L^x_{s_n}\|^{1/2}_{p/2,[a,b]} - A_n,$$

and, by the continuity of $\{L^x_s; 0 \leq s \leq t\}$,

$$(5.41) \qquad \lim_{n \to \infty} \|L^x_{s_n}\|_{p/2,[a,b]} = \|L^x_t\|_{p/2,[a,b]}.$$

This completes the proof of Theorem 1.2 for $-\infty < a < b < \infty$. To handle, for example, $a = -\infty, b = \infty$, note that by what we have shown, almost surely,

$$(5.42) \quad \lim_{h \downarrow 0} \int_{-k}^{k} \frac{|L^{x+h}_t - L^x_t|^p}{\sigma^p_0(h)} \, dx = 2^{p/2} E|\eta|^p \int_{-k}^{k} |L^x_t|^{p/2} \, dx, \qquad k = 1, 2, \ldots.$$

The case $a = -\infty, b = \infty$ follows, since, for each $t$, $L_t^x$ has compact support in $x$ almost surely. $\square$

**Acknowledgments.** This paper was motivated by the work of Richard Bass, Xia Chen and Jay Rosen on large deviations for $L^p$ moduli of continuity of local times.

DEPARTMENT OF MATHEMATICS
CITY COLLEGE, CUNY
NEW YORK, NEW YORK 10031
USA
E-MAIL: mbmarcus@optonline.net

DEPARTMENT OF MATHEMATICS
COLLEGE OF STATEN ISLAND, CUNY
STATEN ISLAND, NEW YORK 10314
USA
E-MAIL: jrosen3@earthlink.net